\documentclass{amsart}

\usepackage{amsmath,amssymb,amsfonts,amsthm,amscd,latexsym,euscript}
\usepackage[all]{xy}

\newcommand{\Hor}{\ensuremath{\textup{Hor}}}
\newcommand{\Rmod}{\ensuremath{R\textup{-mod}}}

\newcommand{\dgrm}[1]{\ensuremath{\smash{\underset{\widetilde{\hphantom{#1}}}{#1}} \mathstrut}}

\newcommand{\mhyphen}{\ensuremath{\mathrel{\mathstrut}}}

\newcommand{\domain}[1]{\ensuremath{\mathrm{dom}({#1})}}
\newcommand{\range}[1]{\ensuremath{\mathrm{codom}({#1})}}

\newcommand{\cal}[1]{\ensuremath{\mathcal #1}}

\newtheorem {theorem1}{Theorem}[section]
\newtheorem {theorem}[theorem1]{Theorem}

\newtheorem {proposition}[theorem1]{Proposition}
\newtheorem {lemma}[theorem1]{Lemma}
\theoremstyle{definition}
\newtheorem {definition}[theorem1]{Definition}
\newtheorem {example}[theorem1]{Example}
\theoremstyle{remark}
\newtheorem {remark}[theorem1]{Remark}

\newcommand{\calT}{\ensuremath{\mathcal{T}}}

\newcommand{\cat}[1]{\ensuremath{\EuScript #1}}

\newcommand{\Top}{\ensuremath{\calT\textup{op}}}

\newcommand{\colim}{\ensuremath{\mathop{\textup{colim}}}}

\def\id{\ensuremath{\textsl{id}}}

\newcommand{\rarrow}{\rightarrow}
\newcommand{\Id}{\ensuremath{\textup{Id}}}

\newcommand{\Vopenka}{Vop\v enka}

\begin{document}

\SelectTips{cm}{10}

%
%

\title[Equivariant cellularization]{Localization with respect to a class of maps II --
                                    Equivariant cellularization and its application}
\author{Boris Chorny}
\address{Einstein Institute of Mathematics, Edmond Safra Campus, Givat
Ram, The Hebrew University of Jerusalem, Jerusalem, 91904, Israel}
\curraddr{Department of Mathematics, Middlesex College, The University of Western Ontario, London, Ontario N6A 5B7, Canada}

\email{bchorny2@uwo.ca}

\thanks{}

\subjclass{Primary 55U35; Secondary 55P91, 18G55} \keywords{model category,
localization, equivariant homotopy}
\date{\today}
\dedicatory{}
\commby{}

\hfuzz=2.5pt

\begin{abstract}
We present an example of a homotopical localization functor which is not a localization with respect to any set of maps. Our example arises from equivariant homotopy theory. The technique of equivariant cellularization is developed and applied to the proof of the main result. \end{abstract}

\maketitle


\section*{Introduction}
Coaugmented idempotent functors, or localizations, occur all over ma\-the\-ma\-tics under different names (e.g., \emph{idempotent monads} \cite{Fakir}, \emph{orthogonal reflections} \cite[1.36]{AR} or \emph{coreflections} on the full subcategory \cite[IV.3]{MacLane}). These are functors $L\colon \cat C\rarrow \cat C$ equipped with a natural transformation $\eta\colon \Id \rarrow L$, such that the natural maps $L_{\eta_X}, \eta_{L_X}\colon LX \rightrightarrows LLX$ are equal isomorphisms. For example, any multiplicative system $S$ in a commutative ring $R$ gives rise to the localization functor $L\colon \Rmod \rarrow \Rmod$ and a natural coaugmentation map $\eta_M\colon M\rarrow L(M)=M\otimes_R R[S^{-1}]$, which justify the name localization. Another example is given by the abelianization functor in the category of groups. 

The purpose of such construction $L$ is to ``forget'' certain information about the mathematical structure one studies. In the case when the information we want to discard may be described in terms of a \emph{set} of generators, the existence of the localization is usually implied by the standard results from the category theory \cite[1.37]{AR}. But if the description of the undesired information is available only in terms of a proper class of generating maps, no general results are available in the standard set-theoretical framework. However, if the underlying category is \emph{locally presentable} \cite[Def.~1.17]{AR}, then it is known that the general question of existence of orthogonal reflection is equivalent to the weak \Vopenka's principle \cite[6.22, 6.23]{AR}. Moreover, if one is ready to assume a more powerful \Vopenka's principle, then by \cite[6.24]{AR} any class of maps (which is closed under limits in the category of maps) may be generated by just a set of maps (i.e., it is a \emph{small-orthogonality} class \cite[Def.~1.32]{AR}). 

\Vopenka's principle and weak \Vopenka's principle are set theoretical statements which are known to be unprovable in the standard set theory, but not known to be independent of the rest of the axioms. See Chapter~6 of J.~Ad\'amek and J.~Rosick\'y's book \cite{AR} for more information.

In this paper we give a counterexample to the similar question about homotopical localization. We do not assume here any of the non-standard axioms. 

Homotopical localizations played an important role in algebraic topology and algebraic geometry over past thirty years. Most of the previously known homotopy idempotent functors (except few rare cases like \cite{Farjoun:towers}, \cite{Isaksen}) are known to be equivalent to the localization with respect to a set of maps, though sometimes it is non-trivial to find an appropriate set (see \cite[1.E]{Farjoun-book} for examples and discussion). This lead E.~Dror Farjoun to ask in \cite{Farjoun:v1}: whether any homotopy idempotent functor on the category of spaces is an $S$-localization for some set of maps $S$? Like in the non-homotopical version, it turned out that this question cannot be settled using only the standard axioms of set theory, but the answer is affirmative if one assumes, additionally, \Vopenka's principle \cite{CSS}. Recently these results were extended to any simplicial combinatorial(= cofibrantly generated \& locally presentable) model category \cite{ChC}.

In this work we show that it is impossible to remove the assumption on the model category to be cofibrantly generated. In more detail, we consider the category of diagrams of spaces with the equivariant model structure \cite{Farjoun}, which is known to be non-cofibrantly generated \cite{Chorny}, but still reasonable enough to admit the standard localization theory with respect to a set of maps and even with respect to a class of maps satisfying some restrictive conditions \cite{PhDI}. The main result of the present paper is that the functor which associates the constant diagram of points to any diagram $\dgrm X$ is not a localization functor with respect to any set of maps of diagrams. This is really a localization with respect to a class of maps.

We would like to stress that our counterexample is not based on some anomaly of the underlying category (in fact, it can be chosen to be locally presentable). Our argument uses the observation that the considered model category is not cofibrantly generated, but it does not mean that our example may be generalized to any non-cofibrantly generated model category: there exists an example of a model category which is Quillen equivalent to the trivial model category, but it is not cofibrantly generated \cite{AHRT}, hence any localization functor in this model category is a localization with respect to an empty set of maps - a fibrant replacement.

Together with the works \cite{CSS, ChC} our example provides an answer to Dror Farjoun's question.

\subsection*{Organization of the paper} 
We continue here the discussion started in \cite{PhDI} and use freely results and notions introduced there. Let us recall only the definition of the equivariant model structure on the category of $D$-shaped diagrams of spaces \cite{Farjoun}: a diagram \dgrm T is called an \emph{orbit} if $\colim_D \dgrm T = \ast$; a map $f\colon \dgrm X \rarrow \dgrm Y$ is a \emph{weak equivalence} or  \emph{fibration} if the induced map $\hom(\dgrm T, f)\colon \hom(\dgrm T, \dgrm X) \rarrow \hom(\dgrm T, \dgrm Y)$ is a weak equivalence or fibration of simplicial sets respectively. Here and further in the paper $\hom(\,\cdot\, , \,\cdot\,)$ denotes the simplicial function complex.

In the fist section we develop the theory of equivariant colocalizations or cellularizations, which is complementary to the theory of equivariant localizations introduced in \cite{PhDI}. We do not treat here the question of the existence of the fixed-pointwise cellularization, since we do not have applications for this notion. The category of spaces may be taken to be the category of simplicial sets or the compactly generated topological spaces with the standard simplicial model structure. 

In the second section we apply the results of the first section in order to prove that the functor $L$ which entirely discards the homotopical information of diagrams of spaces: $L(\dgrm X) = \dgrm\ast$ is not a localization with respect to any set of arrows between diagrams. Our proof works only for the case of diagrams of topological spaces, so that every diagram is fibrant in the equivariant model category. In order to obtain an example of a category which is still locally presentable we may use the category of $I$-generated topological spaces introduced recently by J.~Smith. A useful account of J.~Smith's ideas on the $I$-generated spaces is given by D.~Dugger \cite{Dugger:Delta-Notes}.

\medskip
\paragraph{\emph{Acknowledgements.}} I would like to thank Emmanuel Dror Farjoun for his support and many helpful ideas. I am grateful to Jeff Smith for catching a mistake in an early version of this work.

\section{Construction of the equivariant colocalization functor}\label{coloc}
In this section we explain the complementary approach to the localization.
Namely, we construct for any set of cofibrant diagrams $A = \{\dgrm A\}$ the
augmented functor $CW_{A}\colon \cal S^D\rarrow \cal S^D$ such that for each
$\dgrm X \in \cal S^D, \; CW_{A}(\dgrm X)$ is $A$-colocal and the natural map
$p_{\dgrm X}\colon CW_{A}(\dgrm X)\rarrow \dgrm X$ is an $A$-colocal
equivalence (see below). We prove also that the natural map $p_{\dgrm
X}\colon  CW_{A}(\dgrm X)\rarrow \dgrm X$ is terminal (up to homotopy) among
all the maps of the $A$-colocal spaces into \dgrm X, thus $CW_{A}\dgrm X$ is
characterized up to a weak equivalence.

The crucial difference between the equivariant framework and the ordinary one is that it is not pointless to consider the cellularization of non-pointed diagrams (cf.~\cite[p.~40]{Farjoun-book}, \cite[3.1.10]{Hirschhorn}). Already in the case of a group $G$ acting on topological spaces there are non-trivial augmented homotopy idempotent functors on $\cal S^G$, e.g., for any subgroup $H$ of $G$ there exists the augmented functor which assigns to any $G$-space \dgrm X its subspace fixed by $H$ with the natural inclusion $p_{\dgrm X}\colon (\dgrm X)^H \rarrow \dgrm X$. The homotopy idempotence is clear.

\subsection{Preliminaries on colocal diagrams and colocal equivalences}
\begin{definition}
Let $A$ be a set of cofibrant diagrams.

%

\begin{itemize}
\item
A map $f\colon \dgrm X \rarrow \dgrm Y$  is an $A$-\emph{colocal equivalence}
(or just an $A$-\emph{equivalence}) if for any fibrant replacement $\hat f$ of
$f$ the induced map
\[
\hom(\dgrm A, \hat{f}) \colon \hom(\dgrm A, \dgrm {\hat{X}}) \rarrow \hom(\dgrm
A, \dgrm {\hat{Y}})
\]
is a weak equivalence of simplicial sets for every $\dgrm A\in A$.
\item
A cofibrant diagram $\dgrm B$ is $A$-\emph{colocal} if for any $A$-colocal
equivalence $g\colon \dgrm X\rarrow \dgrm Y$ and any fibrant replacement $\hat
g$ of $g$ the induced map
\[
\hom(\dgrm B, \hat{g})\colon  \hom(\dgrm B, \hat{\dgrm X}) \rarrow \hom(\dgrm
B, \hat{\dgrm Y})
\]
is a weak equivalence of simplicial sets.
\end{itemize}
\end{definition}

\begin{remark}\label{rem}
The above notions are well-defined, i.e., they do not depend on
the choice of the fibrant replacement. It follows from
\cite[9.7.2]{Hirschhorn}. We shall use also an $A$-colocal version
of the Whitehead theorem (see \cite[3.2.13]{Hirschhorn} for the
proof).
\end{remark}

\begin{proposition}[\dgrm A-colocal Whitehead theorem]
A map $g\colon Q_1\rarrow Q_2$ is a weak equivalence of $A$-colocal diagrams if
and only if $g$ is an $A$-colocal equivalence.
\end{proposition}

\begin{proposition}\label{coloc:equiv}
A map $g\colon \dgrm X \rarrow \dgrm Y$ of diagrams is an $A$-colocal
equivalence if and only if there exists a fibrant replacement $\hat g$ of $g$
which has the right lifting property with respect to the following families of
maps:
\begin{itemize}
\item generating trivial cofibrations $J$;
\item
$\Hor(A) = \{\partial\Delta^n\otimes \dgrm A \hookrightarrow \Delta^n
                           \otimes \dgrm A \;|\; n\geq 0,\; \dgrm A\in A\}$.

\end{itemize}
\end{proposition}

\begin{proof}
If $g$ is an $A$-colocal equivalence, then there exists a fibrant replacement
$\hat g\colon \hat{\dgrm X}\rarrow \hat{\dgrm Y}$ such that $\hom( \dgrm A,
\hat g)$ is a weak equivalence for every $\dgrm A\in A$. Consider the
factorization of $\hat g = \hat{g}'i$ into the trivial cofibration $i$ followed
by the fibration $\hat{g}'$, which is also an $A$-equivalence by
Remark~\ref{rem}. Hence, $\hat{g}'$ is a fibrant approximation of $g$, which
has the right lifting property with respect to the elements of $J$, as a
fibration of diagrams, and with respect to the elements of $\Hor(A)$ by
adjunction.

Conversely, if $\hat g$ is a fibrant replacement of $g$ with the right lifting
property with respect to the elements of $J$ and $\Hor(A)$, then, by
adjunction, $\hom(\dgrm A, \hat g)$ is a trivial fibration of simplicial sets
for every $\dgrm A\in A$, therefore $g$ is an $A$-colocal equivalence.
\end{proof}

\begin{proposition}\label{instr:K2}
The class of maps $K = J\cup \Hor(A)$ may be equipped with an instrumentation.
\end{proposition}

The proof of this proposition is similar to the proof of Proposition~
\cite[6.6]{PhDI} and is left to the reader.

\subsection{Construction of $CW_A$}
The naive approach which is dual to the construction of the localization does
not work: if for any diagram \dgrm X we factor the map $\emptyset\rarrow \dgrm
X$ into a $K$-cellular map followed by a $K$-injective map, then for
non-fibrant \dgrm X we will not be able to show that the $K$-injective map is
an $A$-equivalence. See \cite[5.2.7]{Hirschhorn} for a counterexample in the
category of pointed simplicial sets.

The right properness of the category of diagrams is essential for
the following construction (compare \cite[5.3.5]{Hirschhorn}). For
any diagram \dgrm X choose a functorial cofibrant fibrant
approximation $\smash{j\colon \dgrm X
\mathbin{\tilde{\hookrightarrow}} \hat{\dgrm X}}$. Apply the
generalized small object argument, with respect to the
instrumented (by Proposition~\ref {instr:K2}) class $K$, to
factorize the map $\smash{\emptyset \rarrow \hat{\dgrm X}}$ into a
$K$-cellular map $r$ followed by a $K$-injective map $s$:
\[
\emptyset \overset r \longrightarrow \hat{\dgrm W} \overset s \longrightarrow
\dgrm X.
\]
Next, take $\dgrm W = \dgrm X \times_{\hat{\dgrm X}} \hat{\dgrm W}$; then the
natural map $t\colon \dgrm W \rarrow \hat{\dgrm W}$ is a weak equivalence as a
pullback of a weak equivalence $j$ along the fibration $s$ in the right proper
model category of diagrams. ($s$ is a fibration, since $s\in K$-inj.) The
natural map $v\colon \dgrm W\rarrow \dgrm X$ is in $K$-inj as a pullback of
the $K$-injective map $s$. The functorial fibrant cofibrant approximation
$\emptyset \hookrightarrow CW_{A}\dgrm X \mathbin{\underset u
{\tilde{\twoheadrightarrow}}} \dgrm W$ supplies us with an augmented functor
$CW_{A}\dgrm X$, where the augmentation $p_{\dgrm X}$ is given by the
composition
\[
CW_{A}\dgrm X \xrightarrow{u} \dgrm W \xrightarrow{v} \dgrm X, \; p_{\dgrm
X}=vu.
\]

Summarizing, we have the commutative diagram
\[
\xymatrix{
 \emptyset \ar@{^{(}->}@(d,l)[dr]_g \ar@{^{(}->} @(r, ul) [drrr]^r
 \\
    & CW_{A} \dgrm X \ar@{->>}[r]^\sim_u \ar@{->>}[dr]_{p_{\dgrm X}}& \dgrm W \ar[r]^\sim_t \ar@{->>}[d]^v & {\hat{\dgrm W}} \ar@{->>}[d]^s \\
    &                       & \dgrm X \ar@{^{(}->}[r]^\sim_j & {\hat{\dgrm X}}\ar@{->>}[r] & \ast\ldotp\\
}
\]

The map $p_{\dgrm X}\colon CW_{A}\dgrm X \rarrow \dgrm X$ is an $A$-colocal
equivalence by Proposition~\ref{coloc:equiv}, since its fibrant approximation
$s\colon \hat{\dgrm W}\rarrow \hat{\dgrm X}$ is $K$-injective, i.e., it has the
right lifting property with respect to the sets $J$ and $\Hor(A)$.
\par
\begin{remark}\label{rem:K-inj}
We note, for future reference, that $p_{\dgrm X}\in K$-inj. $p_{\dgrm X} = vu$
is a composition of two fibrations, hence a fibration, i.e., it has the right
lifting property with respect to any element of $J$. For any element $\dgrm
C\hookrightarrow \dgrm D$ of $\Hor(A)$ and any commutative square
\[
\xymatrix{
 {\dgrm C}
 \ar[r]
 \ar@{^{(}->}[d]
                  &  {CW_{A}\dgrm X}
                        \ar@{->>}[d]^{p_{\dgrm X}}\\
 {\dgrm D}
 \ar[r]           &  {\dgrm X}\\
}
\]
we construct first a lift $\hat{h}\colon \dgrm D\rarrow \hat{\dgrm W}$, which
exists since $s\in K$-inj. Let $h\colon \dgrm D \rarrow \dgrm W$ be the
natural map into the pullback \dgrm W. Finally, the required lift $l\colon
\dgrm D \rarrow CW_{A}\dgrm X$ exists, since the map $u$ is a trivial
fibration and any element of $\Hor(A)$ is a cofibration.
\end{remark}

It remains to show that $CW_{A}\dgrm X$ is $A$-colocal for any diagram \dgrm X.
But $CW_{A}\dgrm X$ is cofibrant and weakly equivalent to the $K$-cellular
complex $\hat{\dgrm W}$, hence it will suffice to show that $\hat{\dgrm W}$ is
$A$-colocal. But any $K$-cellular diagram is $A$-colocal. The following
proposition completes the proof.
\par
\begin{proposition}
Any $K$-cellular complex \dgrm B is an $A$-colocal diagram.
\end{proposition}
\begin{proof}
We will prove this by the transfinite induction on the indexing ordinal of the
$\lambda$-sequence $\emptyset=\dgrm B_0\hookrightarrow \dgrm B_1\hookrightarrow
\cdots\hookrightarrow \dgrm B_i \hookrightarrow \cdots$, whose colimit is \dgrm
B.

The diagram $\emptyset$ is obviously $A$-colocal for any set of diagrams $A$,
hence the base of the induction.

For each $k$ the map $i_k\colon \dgrm B_k \hookrightarrow \dgrm
B_{k+1}$ is a pushout of a coproduct of the elements of $K$. Each
map in $K$ has the homotopy left lifting property with respect to
the class of fibrations which are $A$-equivalences. For the
elements of $\Hor(A)$ this follows from
\cite[9.4.8(1)]{Hirschhorn} (substitute $i$ by $\emptyset
\rightarrow \dgrm A$, $\dgrm A \in A$ and $(K,L) =
(\partial\Delta^n, \Delta^n)$). Maps of $J$ are trivial
cofibrations, therefore they have the homotopy left lifting
property with respect to all fibrations \cite[9.4.4]{Hirschhorn}.
Next, a coproduct of a set of maps with homotopy left lifting
property with respect to any class of maps again has the homotopy
left lifting property with respect to the same class: first note
that coproducts commute with pushouts and with the left adjoint
functors $\cdot \otimes K$ for any simplicial set $K$, then apply
\cite[9.4.7(2)]{Hirschhorn}. But the homotopy left lifting
property is preserved under pushouts, hence the map $i_k\colon
\dgrm B_k \hookrightarrow \dgrm B_{k+1}$ has the homotopy left
lifting property with respect to any fibration which is an
$A$-equivalence.

First we prove the inductive step for successor ordinals. Suppose $\dgrm B_k$
is \dgrm A-colocal. For any \dgrm A-equivalence $g\colon \dgrm X\rarrow \dgrm
Y$ take $\hat{g}\colon \hat{\dgrm X}\rightarrow \hat{\dgrm Y}$ to be its
fibrant approximation. In the commutative diagram
\[
\xymatrix{
  {\hom(\dgrm B_{k+1}, \hat{\dgrm X})}
  \ar@{->>}[dd]_{\hom(i_k, \hat{\dgrm X})}
  \ar[rr]^\sim_{\hom(\dgrm B_{k+1}, \hat g)}
  \ar[dr]^{\dir{~}}_{m}
                                       &           & {\hom(\dgrm B_{k+1}, \hat{\dgrm Y})}
                                                      \ar@{->>}[dd]^{\hom(i_k, \hat{\dgrm Y})}\\
                                       & {\dgrm P}
                                          \ar[ur]^{\dir{~}}_n
                                          \ar@{->>}[dl]\\
  {\hom(\dgrm B_{k}, \hat{\dgrm X})}
   \ar[rr]^\sim_{\hom(\dgrm B_k, \hat{g})}
                                       &           & {\hom(\dgrm B_{k},   \hat{\dgrm Y})}\\
 }
\]
the map $\hom(\dgrm B_k, \hat{g})$ is a weak equivalence of simplicial sets by
inductive assumption. The map $\hom(i_k, \hat{\dgrm Y})$ is a fibration, since
$i_k$ is a cofibration and $\hat{\dgrm Y}$ is a fibrant diagram. Let
\[
\dgrm P = \hom(\dgrm B_{k}, \hat{\dgrm X}) \times_ {\hom(\dgrm B_{k},
\hat{\dgrm Y})} \hom(\dgrm B_{k+1}, \hat{\dgrm Y});
\]
therefore the map $n$ is a weak equivalence as a pullback of the weak
equivalence $\hom(\dgrm B_k, \hat{g})$ along the fibration $\hom(i_k, \hat{
\dgrm Y})$. The map $m$ is a weak equivalence, since $i_k$ has the homotopy
left lifting property with respect to the map $\hat g$, which is a fibration
and an $A$-equivalence. Finally, $\hom(\dgrm B_{k+1}, \hat g) = n\circ m$ is a
weak equivalence as a composition of weak equivalences $m$ and $n$. This
proves the inductive step for successor ordinals.

If $\mu \leq \lambda$ is a limit ordinal, then $\dgrm B_\mu = \colim_{k<\mu}
\dgrm B_k$. For any $A$-equivalence $g\colon \dgrm X\rarrow \dgrm Y$, let $\hat
g\colon \hat{\dgrm X} \rarrow \hat{\dgrm Y}$ be its fibrant approximation. If
$\dgrm B_k$ is $A$-colocal for each $k < \mu$, then the induced map of
simplicial sets $\hom(\dgrm B_k, \hat g)$ is a weak equivalence. This implies
that $\hom(\dgrm B_\mu, \hat g) = \hom(\colim_{k < \mu}\dgrm B_k, \hat g) =
\lim_{k < \mu}\hom(\dgrm B_k, \hat g)$ is a weak equivalence of simplicial
sets, since $\lim_{k < \mu}\hom(\dgrm B_k, \hat{\dgrm X})$ and $\lim_{k <
\mu}\hom(\dgrm B_k, \hat{\dgrm Y})$ are homotopy inverse limits of towers of
fibrations.

Hence, the step of the induction.
\end{proof}

\subsection{Universality of $CW_A$}
%
%
\begin{proposition}[$CW_{A}$ is terminal]
For any map $u\colon \dgrm U\rarrow \dgrm X$ from an $A$-colocal diagram \dgrm
U there exists a factorization $\dgrm U\rarrow CW_{A}(\dgrm X)\rarrow \dgrm X$
which is unique up to simplicial homotopy.
\end{proposition}
\begin{proof}
By Remark \ref{rem:K-inj} the natural map $p_{\dgrm X}\colon
CW_{A}\dgrm X \rarrow \dgrm X$ is in $K$-inj. The map $\emptyset
\rarrow \dgrm U$ is in $K$-cof, since \dgrm U is $A$-colocal (one
adopts the proof of \cite[3.4.1]{Hirschhorn} in our case). Then
the following commutative square
\[
\xymatrix{
 {\emptyset}
 \ar[r]
 \ar@{^{(}->}[d]
                  &  {CW_{A}\dgrm X}
                        \ar@{->>}[d]^{p_{\dgrm X}}\\
 {\dgrm U}
 \ar[r]_q
 \ar@{-->}[ur]
                  &  {\dgrm X}\\
}
\]
admits a lift, which provides the required factorization.

The factorization above is unique, since the map $\hom(\dgrm U,
p_{\dgrm X })$ is a weak equivalence by
\cite[13.2.2(2)]{Hirschhorn}, thus the induced map on simplicial
homotopy classes is a bijection and, in particular, injection,
therefore non-homotopic lifts cannot correspond to the same class
of maps $[q] \in [\dgrm U, \dgrm X]$.
\end{proof}

Like the localization functor, $CW_{A}$ has also the second universal
property, more precisely, its restriction to the subcategory of fibrant
diagrams does. The augmentation map $p_{\dgrm X}$ is a fibration for any \dgrm
X, hence the subcategory of fibrant diagrams is stable under localizations.
Denote by $CW_{A}^r$ the restriction of $CW_{A}$ on the subcategory of fibrant
diagrams (do nothing for the diagrams of topological spaces). Then $CW_{A}^r$
is initial with respect to the $A$-equivalences. In more detail, we have the
following
\begin{proposition}[$CW_{A}^r$ is initial]
On the subcategory of fibrant diagrams the augmentation map $p_{\dgrm X}\colon
CW_{A}^r(\dgrm X)\rarrow \dgrm X$ is initial, up to homotopy, among all
$A$-colocal equivalences, i.e., for any $A$-equivalence of fibrant diagrams
$f\colon \dgrm Y\rarrow \dgrm X$, there exists a unique, up to homotopy, map
$g\colon CW_{A}^r(\dgrm X)\rarrow \dgrm Y$ such that $p_{\dgrm X} \overset s
\sim fg$.
\end{proposition}
\begin{proof}
Apply the functor $CW_{A}^r$ on the $A$-colocal equivalence $f$,
then the map $CW_{A}^r(f)=CW_{A}(f)$ is an $A$-equivalence by the
`2 out of 3' property for $A$-equivalences
\cite[3.2.3]{Hirschhorn}. The $A$-colocal Whitehead theorem
implies that $CW_{A}(f)$ is a weak equivalence, therefore it has a
homotopy inverse $q$. If we take  $g = p_{\dgrm Y}q$, then $fg =
fp_{\dgrm Y}q = p_{\dgrm X}CW_{A}(f)q \overset s \sim p_{\dgrm
X}\id_{CW_{A}(\dgrm X)} = p_{\dgrm X}$.
\[
\xymatrix{
 {CW_{A}\dgrm Y}
  \ar[r]^{CW_{A}(f)}_\sim
  \ar[d]_{p_{\dgrm Y}}
                            &  {CW_{A}\dgrm X}
                              \ar[d]_{p_{\dgrm X}}
                              \ar@{-->}[dl]^g
                              \ar@{.>}@/_20pt/[l]|q\\
 \dgrm Y
  \ar[r]_f
                            &  \dgrm X}
\]

Suppose there exists $g'\colon CW_{A}(\dgrm X)\rarrow \dgrm Y, \;
g'\neq g$, and such that $fg' \overset s \sim p_{\dgrm X}$. By the
terminal property of the cellularization functor there exists a
map $q'\colon CW_{A}(\dgrm X) \rarrow CW_{A}(\dgrm Y)$ such that
$g' = p_{\dgrm Y}q'$. It will suffice to show that $q' \overset s
\sim q$, since it implies that $g' \overset s \sim g$ by
\cite[9.5.4]{Hirschhorn}. The induced map on simplicial homotopy
classes $\psi = [CW_{A}(\dgrm X), fp_{\dgrm Y}]$ is a bijection,
since $fp_{\dgrm Y}$ is an $A$-equivalence of fibrant diagrams and
$CW_{A}(\dgrm X)$ is $A$-colocal. But $\psi([q]) = \psi([q'])$,
because $fp_{\dgrm Y}q = fg \overset s \sim p_{\dgrm X} \overset s
\sim fg' = fp_{\dgrm Y}q'$, hence $[q] = [q']$ or $q' \overset s
\sim q$.
\end{proof}

\section{Application: An example of a coaugmented, homotopy idempotent functor
         that is not a localization with respect to a set of maps}\label{example}
Assuming \Vopenka's principle, every continuous localization functor in a combinatorial, simplicial model category is a localization with respect to a set of maps \cite{ChC}. May this property be generalized? No, in this section we introduce a counterexample.

We will prove that if the index category $D$ has the proper class of orbits
(which satisfy some mild conditions, see below), then the simplest
localization functor which associates to any diagram \dgrm X the final object
$\ast = \textsl{pt}$ is not equivalent to an $S$-localization functor for any
set of maps $S$ of diagrams, i.e., there is no set of maps $S$ such that
$L_S(\dgrm X)$ is contractible for every diagram \dgrm X. The current proof
works only for $\cal S$ being the category of compactly generated spaces or $I$-generated spaces \cite{Dugger:Delta-Notes}, since we use the assumption that all the objects of $\cal S^D$ are fibrant.

This example is a continuation of Example~\cite[7.1]{PhDI}, where the
fixed-pointwise localization of diagrams of spaces with respect to the map
$f\colon \emptyset \rarrow \ast$ of spaces was considered. This is essentially
the localization with respect to the class $F$ of inclusions of the empty
diagram into the orbits of $\cat O$. The class of horns on $F$ in our case is
$\Hor(F) = \{\partial \Delta^n \otimes \dgrm T \hookrightarrow \Delta^n
\otimes \dgrm T \;|\; n\geq 0,\; \dgrm T\in \cat O\}$. Note that $\Hor(F) = I$,
the class of generating cofibrations. This suggests that the proof would be
similar to the proof that the model category on the category of diagrams of
spaces is not cofibrantly generated \cite{Chorny}. In fact, the method
developed in \cite{Chorny} works: we find a class of orbits each of which is a
retract of a diagram built out of a fixed \underline{set} of cofibrations, and
this is a contradiction to \cite[2.2]{Chorny}.

\begin{definition}
An orbit \dgrm T is \emph{regular} if there exists a map $\ast \rarrow \dgrm
T$ or, equivalently, $\lim \dgrm T \neq \emptyset$. Otherwise, \dgrm T is
\emph{singular}.
\end{definition}

\begin{example}
If $D = J = (\bullet \rarrow \bullet)$ is the category with two objects and one
non-identity morphism, then the only singular orbit is $\dgrm T_0 = (\emptyset
\rarrow \ast)$. The rest of the orbits are regular. In this case we have a
proper class of regular orbits; this is precisely the condition on the
category $D$ that guarantees that the functor $L(\dgrm X) = \ast$ is not an
$S$-localization for any set of maps $S$.
\end{example}

\begin{example}
If $D = (\bullet \qquad \bullet \rarrow \bullet)$ is the category with three
objects and one non-identity morphism, then there are no regular orbits at
all, though there is a proper class of orbits $\{(\ast \qquad \emptyset
\rarrow \emptyset);\; (\emptyset \qquad X\rarrow \ast) \;|\; X\in \cal S\}$.
For this indexing category we do not know whether the functor $L(\dgrm X) =
\ast$ is equivalent to some $L_S$.
\end{example}

\begin{theorem}\label{main-theorem}
Let $D$ be a small category that satisfies the following condition: the
category $\cat O_D$ of $D$-orbits contains a proper class of regular orbits.
Then there does not exist a set of maps $S$ of $D$-diagrams such that
$L_S(\dgrm X) \simeq \ast$ for every diagram $\dgrm X$.
\end{theorem}
\begin{lemma}\label{lemma-one}
Suppose that there exists a set $B$ of cofibrant diagrams with the following
property: a fibrant diagram \dgrm Z is contractible if and only if the
simplicial set $\hom(\dgrm B, \dgrm Z)$ is contractible for every $\dgrm B\in
B$. Then every regular $D$-orbit \dgrm T is $B$-colocal.
\end{lemma}
\begin{proof}
Every orbit \dgrm T is cofibrant in the model category generated by orbits,
hence it is enough to show that if \dgrm T is regular, then for any $B$-colocal
equivalence $g\colon \dgrm X \rarrow \dgrm Y$ and for any fibrant replacement
$\hat g\colon \hat{\dgrm X} \rarrow \hat{\dgrm Y}$ of $g$, the induced map
\[
\hom(\dgrm T, \hat g)\colon \hom(\dgrm T, \hat{\dgrm X})\rarrow \hom(\dgrm T,
\hat{\dgrm Y})
\]
is a weak equivalence of simplicial sets.

Without loss of generality we may assume that $\hat g$ is a fibration. Take
any map $\ast \rarrow \hat{\dgrm Y}$. It may happen that such a map does not
exist, in which case there are no maps $\dgrm T \rarrow \hat{\dgrm Y}$, since
\dgrm T is a regular orbit and admits a map $\ast \rarrow \dgrm T$. There are
also no maps $\dgrm T \rarrow \hat{\dgrm X}$, otherwise it could be
concatenated with $\hat g$. In this case we are done, since $\hom(\dgrm T,
\hat g)$ is the identity map of $\emptyset$, i.e., a weak equivalence.

If the map $\ast \rarrow \hat{\dgrm Y}$ exists, then let $\dgrm F = \ast
\times_{\hat{\dgrm Y}}\hat{\dgrm X}$. According to the terminology of
\cite{Hirschhorn}, $\dgrm F$ is the homotopy fiber of $\hat g$ over the point
$\ast \rarrow \hat{\dgrm Y}$, since the pullback square
\[
\xymatrix{
       {\dgrm F}
        \ar[r]
        \ar@{->>}[d]_h
                  &  {\hat{\dgrm X}}
                        \ar@{->>}[d]^{\hat g}\\
       {\ast}
        \ar[r]
                  &  {\hat{\dgrm Y}}\\
}
\]
is a homotopy fiber square by \cite[13.3.8]{Hirschhorn}.

For every $\dgrm B\in B$ we apply the functor $\hom(\dgrm B, \cdot)$ on the
commutative square of fibrant diagrams above and obtain the pullback of
simplicial sets:
\[
\xymatrix{
       & {\hom(\dgrm B, {\dgrm F})}
          \ar[r]
          \ar[d]_{\hom(\dgrm B, h)}
                  &  {\hom(\dgrm B, \hat{\dgrm X})}
                        \ar@{->>}[d]^{\hom(\dgrm B, \hat g)}\\
 {\ast}
 \ar@{=}[r]^>(.6){\sim}
     & {\hom(\dgrm B, \ast)}
       \ar@{->>}[r]
                  &  {\hom(\dgrm B, \hat{\dgrm Y}).}\\
}
\]

By the assumption, $g$ is a $B$-colocal equivalence, therefore $\hom (\dgrm B,
\hat g)$ is a weak equivalence for every $B\in \dgrm B$, hence $\hom(\dgrm B,
h)$ is a weak equivalence of simplicial sets, since trivial fibrations are
preserved by pullbacks. In other words, $\hom(\dgrm B, {\dgrm F})$ is
contractible for every $\dgrm B\in B$, but the assumption on the set $B$
implies that ${\dgrm F}$ is contractible.

Next, apply the functor $\hom(\dgrm T, \cdot)$ on the initial pullback square.
We obtain a pullback square of simplicial sets, which is also a homotopy fiber
square
\[
\xymatrix{
        {\hom(\dgrm T, \dgrm F)}
         \ar[r]
         \ar[d]
                  &  {\hom(\dgrm T, \hat{\dgrm X})}
                        \ar@{->>}[d]^{\hom(\dgrm T, \hat g)}\\
        {\ast}
         \ar[r]^x
                  &  {\hom(\dgrm T, \hat{\dgrm Y}).}\\
}
\]
The map $\hom(\dgrm T, \dgrm F)\rarrow \ast$ is a weak equivalence, since
$\dgrm F$ is contractible, hence the homotopy fiber of the map $\hom(\dgrm T,
\hat g)$ over the point $x$ is contractible.

We obtain the same conclusion if we consider the homotopy fibre of $\hom(\dgrm
T, \hat{g})$ over any connected component of $\hom(\dgrm T, \hat{\dgrm Y})$.
The point $x$ may be chosen in any component, since every $0$-simplex of
$\hom(\dgrm T, \hat{\dgrm Y})$ corresponds to a map $\dgrm T \rarrow \hat
{\dgrm Y}$; but $\dgrm T$ is a regular orbit, therefore we can choose $x =
\hom(\dgrm T, k)$, where $k$ is the composition $\ast \rarrow \dgrm T \rarrow
\hat{\dgrm Y}$. But a map of simplicial sets with the contractible homotopy
fiber over each component of the base is a weak equivalence! Since $g$ was an
arbitrarily chosen $B$-colocal equivalence, \dgrm T is $B$-colocal.
\end{proof}

In the next lemma we were unable to get rid of the assumption $\cal S = \Top$.
One should expect that it could be replaced by the right properness.
\begin{lemma}\label{lemma-two}
Let $B$ be a set of cofibrant diagrams in the model category of diagrams of
topological spaces generated by the collection of orbits. If $\Hor(B) =
\{\partial \Delta^n \otimes \dgrm B \hookrightarrow \Delta^n \otimes \dgrm B
\;|\; n\geq 0, \dgrm B\in B\}$, then:
\begin{enumerate}
\item
Any map $g\colon \dgrm X \rarrow \dgrm Y$ in $\Hor(B)$-inj is a $B$-colocal
equivalence;
\item
Any $B$-colocal diagram \dgrm X is a retract of a $\Hor(B)$-cellular diagram
$C(\dgrm X)$.
\end{enumerate}
\end{lemma}
\begin{proof}
The advantage of the diagrams of topological spaces is that every object is
fibrant, hence the fibrant replacement functor may be chosen to be the
identity. The first claim follows by applying adjunction and concluding that
for every $\dgrm B\in B$, $\hom(\dgrm B, g)\colon \hom(\dgrm B, \dgrm X)
\rarrow \hom(\dgrm B, \dgrm Y)$ is a trivial fibration, hence a weak
equivalence.

The second claim follows from the following commutative diagram:
\[
\xymatrix{
        {\emptyset}
         \ar@{^(->}[rr]^s
         \ar@{^(->}[dd]
         \ar@{^(->}[dr]^k
                    &        &  {C(\dgrm X)}
                                \ar@{->>}[dd]^t\\
                     &{CW_{\dgrm B}(\dgrm X)}
                       \ar@{..>}[ur]
                       \ar@{->>}[dr]^{p_{\dgrm X}}\\
        {\dgrm X}
         \ar@{-->}[ur]
         \ar[rr]_{\id_{\dgrm X}}
                     &       &  {\dgrm X \ldotp}\\
}
\]
In this diagram the map $s$ is in $\Hor(B)$-cell and $t$ is in $\Hor(B)$-inj,
for they are obtained by the application of the small object argument on the
map $\emptyset \hookrightarrow \dgrm X$ with respect to the set of maps
$\Hor(B)$. The maps $k$ and $p_{\dgrm X}$ are obtained upon application of the
$CW_{B}(\;\cdot\;)$ functor on \dgrm X.

The dashed arrow exists by the terminal property of the $CW_{B}(\;\cdot\;)$
functor, since the diagram \dgrm X is $B$-colocal. The dotted arrow exists by
the initial property of the $CW_{B}(\;\cdot\;)$ functor, since the map $t\in
\Hor(B)$-inj and, hence, a $B$-colocal equivalence by the first claim (recall
that all the diagrams are fibrant). Let us denote the composition of the dashed
with the dotted arrows by $i$.

Summarizing, \dgrm X is a retract of the $\Hor(B)$-cellular diagram $C(\dgrm
X)$, since the composition $\dgrm X \xrightarrow{i} C(\dgrm X) \xrightarrow{t}
\dgrm X$ is the identity on \dgrm X.
\end{proof}

Finally, we are able to prove our main result.

\begin{proof}[Proof of Theorem \ref{main-theorem}]
Suppose that there exists a set of maps (which are, without loss of
generality, cofibrations between cofibrant diagrams) $S=\{f\colon \dgrm A
\hookrightarrow \dgrm B\}$ such that the localization functor $L_S\colon \cal
S^D \rarrow \cal S^D$ associates a contractible space to each diagram. This
means that $\dgrm X$ is $S$-local if and only if \dgrm X is contractible.

Consider the set of (cofibrant) diagrams $B=\{\dgrm A, \dgrm B \;|\; \dgrm A =
\domain{f}, \dgrm B = \range{f},\mhyphen f\in S\}$. The set $B$ has the
following property: a diagram \dgrm X is contractible if and only if
$\hom(\dgrm B, \dgrm X)$ is contractible for all $\dgrm B\in B$. The `only if'
direction being clear, the inverse direction follows from the fact that such
\dgrm X must be $S$-local, since any map between contractible simplicial sets
is a weak equivalence of spaces.

Hence, the set $B$ of cofibrant diagrams satisfies the assumptions of Lemma~
\ref{lemma-one}. Then every regular orbit \dgrm T is $B$-colocal (recall that
there is a proper class of regular orbits).

By Lemma \ref{lemma-two}, every regular orbit is a retract of some
$\Hor(B)$-cellular diagram. But $\Hor(B)$ is a \underline{set} of cofibrations,
therefore, by \cite[Lemma 2.2]{Chorny}, absolute $\Hor(B)$-cellular complexes
contain only a set of orbits, hence the contradiction.
\end{proof}

\bibliographystyle{abbrv}
\bibliography{xbib}

\end{document}